\documentclass{IOS-Book-Article}

\usepackage{mathptmx}

\usepackage{amsfonts,amsmath,url, amsthm}

\newtheorem{theorem}{Theorem}[section]
\newtheorem{conjecture}[theorem]{Conjecture}

\newcommand{\CC}{\mathbb{C}}
\newcommand{\FF}{\mathbb{F}}
\newcommand{\QQ}{\mathbb{Q}}
\newcommand{\RR}{\mathbb{R}}
\newcommand{\ZZ}{\mathbb{Z}}
\newcommand{\frako}{\mathfrak{o}}
\newcommand{\frakq}{\mathfrak{q}}

\DeclareMathOperator{\Alt}{A}
\DeclareMathOperator{\Aut}{Aut}
\DeclareMathOperator{\Conj}{Conj}
\DeclareMathOperator{\Cyc}{C}
\DeclareMathOperator{\Dih}{D}
\DeclareMathOperator{\End}{End}
\DeclareMathOperator{\Gal}{Gal}
\DeclareMathOperator{\GL}{GL}
\DeclareMathOperator{\SO}{SO}
\DeclareMathOperator{\Sp}{Sp}
\DeclareMathOperator{\ST}{ST}
\DeclareMathOperator{\SU}{SU}
\DeclareMathOperator{\Sym}{S}
\DeclareMathOperator{\USp}{USp}

\begin{document}
\begin{frontmatter}              

\title{Sato-Tate groups of genus 2 curves}
\runningtitle{Sato-Tate groups}

\author{\fnms{Kiran S.} \snm{Kedlaya}%
\thanks{Besides support to participate in the NATO ASI, the author received financial support from NSF grant DMS-1101343 and UC San Diego
(Stefan E. Warschawski Professorship), and was in residence at MSRI during fall 2014 (supported by NSF grant 0932078 000).  Thanks to Jean-Pierre Serre for feedback on a draft.}}

\runningauthor{Kiran S. Kedlaya}
\address{Department of Mathematics, University of California, San Diego, La Jolla, CA 92093, USA}

\begin{abstract}
We describe the analogue of the Sato-Tate conjecture for an abelian variety over a number field; this predicts that the zeta functions of the reductions over various finite fields, when properly normalized, have a limiting distribution predicted by a certain group-theoretic construction related to Hodge theory, Galois images, and endomorphisms. 
After making precise the definition of the \emph{Sato-Tate group} appearing in this conjecture, we describe the classification of Sato-Tate groups of abelian surfaces due to Fit\'e--Kedlaya--Rotger--Sutherland.
\end{abstract}

\begin{keyword}
Sato-Tate group, abelian varieties, equidistribution, Frobenius eigenvalues.
\end{keyword}
\end{frontmatter}

\thispagestyle{empty}
\pagestyle{empty}

\section*{Introduction}
These lecture notes provide an introduction to the theory of \emph{Sato-Tate groups} associated to abelian varieties, with an emphasis on abelian surfaces (especially the Jacobians of hyperelliptic genus 2 curves).

Lecture 1 consists of an introduction to the Sato-Tate conjecture, including some detailed discussion of equidistribution in compact Lie groups. In the manner of \cite{serre-lectures} (which we recommend highly as a gentle introduction for anyone encountering this circle of ideas for the first time), we use the Chebotarev density theorem as a motivating analogy.

Lecture 2 consists of a detailed description of Serre's generalization of the Sato-Tate conjecture to an arbitrary abelian variety \cite{serre-motivic}, as further explicated in \cite{serre-lectures} 
and \cite{banaszak-kedlaya} in terms of the construction of a \emph{Sato-Tate group} using Hodge theory and $\ell$-adic Tate modules. We also give a sketch of Serre's reduction of the generalized Sato-Tate conjecture to some analytic properties of $L$-functions \cite[Chapter~1, Appendix]{serre-abelian} and an indication of how this reduction is used in some known cases of the conjecture (including the original Sato-Tate conjecture over $\QQ$).

Lecture 3 is a survey of the classification of Sato-Tate groups associated to abelian surfaces, due to Fit\'e, Rotger, Sutherland, and the author \cite{fkrs}, including some explicit discussion about the connected parts and component groups which arise. We end with some speculation about generalizations.

\section{The Sato-Tate conjecture}

In this lecture, we recall the statement of the Sato-Tate conjecture (now a theorem!) for elliptic curves, then reformulate it in a fashion which is more suggestive of generalizations.

\subsection{The Hasse bound}

We start with Hasse's theorem on the number of points on an elliptic curve over a finite field.
(see for instance \cite[Theorem~V.1.1]{silverman}). 
\begin{theorem}[Hasse]
Let $E$ be an elliptic curve over a finite field $\FF_q$. Then
\[
\#E(\FF_q) = q+ 1 - a_q, \qquad \left| a_q \right| \leq 2 \sqrt{q}.
\]
\end{theorem}

For example, let $E$ be the projective curve defined by the affine Weierstrass equation
\[
y^2 = x^3 + Ax + B;
\]
then $E$ is an elliptic curve provided that the discriminant $\Delta = -16 (4A^3+27B^2)$ is nonzero. (For $q$ a power of a prime $p \geq 5$, every elliptic curve over $\FF_q$ can be put in this form; see \cite[Chapter~3]{silverman}.) In this case, the set
$E(\FF_p)$ contains exactly one point at infinity (namely the point $[0:1:0]$), so the number of solutions of the equation $y^2 = x^3 +Ax+B$ is $q-a_q$. Hasse's theorem may then be interpreted as follows. 
For an individual value of $x \in \FF_q$, the number of square roots of $x^3 + Ax + B$ can be either 0 or 2, each with probability about 1/2 (there is also a small chance of the number being 1, which we neglect). If the counts were truly independent random variables for distinct values of $x$, the central limit theorem would assert that $\#E(\FF_q) = q+1 - a_q$ where $a_q$ is with high probability less than a constant multiple of $\sqrt{q}$.
Hasse's theorem says that this in fact occurs not just with high probability, but with absolute certainty!

\subsection{Warmup: fixed prime, varying curve}

It is natural to ask how the quantity $a_q$ varies within the interval $[-2\sqrt{q}, 2\sqrt{q}]$ for a fixed $q$. 
Although this is not the question we will be pursuing later, it provides an opportunity to introduce some key ideas in a somewhat simpler setting, so the detour will be worthwhile.

To simplify matters, let us assume $q = p$ is a prime and also that $p \geq 5$, so that we can describe all elliptic curves using Weierstrass equations.\footnote{To do this properly for general $q$, one should work over the set of isomorphism class of elliptic curves modulo $p$, where each class carries the weight $1/w$ for $w$ equal to the number of automorphisms of an elliptic curve in the class. This type of weighting is consistent with the Burnside-Cauchy-Frobenius counting lemma for orbits of group actions.}
Let us view $a_p$ as a random variable over the set of pairs $(A,B) \in \FF_p \times \FF_p$ for which $\Delta \neq 0$,
and then ask questions about the distribution of this random variable. As usual in probability theory, it is useful to study the \emph{moments}
\[
M_d(a_p) = \mathbb{E}(a_p^d) \qquad (d=1,2,\dots)
\]
of $a_p$ (where $\mathbb{E}$ denotes the expected value, i.e., the average over the sample space). Note that $M_d(a_p) = 0$ whenever $d$ is odd, because each elliptic curve can be paired with its quadratic twist to obtain complete cancellation.

The even moments of $a_p$ were studied by Birch \cite{birch}. For small $d$ (and $p \geq 5$), Birch obtains explicit formulas\footnote{The formulas in \cite{birch} include the cases where $\Delta=0$, which we have subtracted out here. Moreover, Serre points out that the formulas in \cite[Theorem~2]{birch} are all missing a factor of $p-1$.} which are polynomials in $p$.
\begin{theorem}[Birch] 
We have
\begin{align*}
M_2(a_p) &= p - p^{-1} \\
M_4(a_p) &= 2p^2 - 3 - p^{-1} \\
M_6(a_p) &= 5p^3 - 9p - 5 - p^{-1} \\
M_8(a_p) &= 14p^4 - 28p^2 - 20p - 7-p^{-1}.
\end{align*}
\end{theorem}
For larger $d$, similar formulas always exist (by the Selberg trace formula) but will not be polynomials in $p$; there are additional contributions coming from the coefficients of certain modular forms. Birch provides one explicit example.
\begin{theorem}[Birch] 
Let $\tau$ denote Ramanujan's $\tau$-function, so that
\[
\sum_{n=1}^\infty \tau(n) q^n = q \prod_{n=1}^\infty (1-q^n)^{24}.
\]
Then
\[
M_{10}(a_p) = 42p^5 - 90p^3 -75p^2 - 35p - 9 - (1 + \tau(p)) p^{-1}.
\]
\end{theorem}

Even when exact formulas for the moments prove to be complicated, it is also useful
to give asymptotic formulas; this will give information about the limiting distribution of $a_p$ as $p \to \infty$ (more on this below).
In this vein, Birch shows the following.
\begin{theorem}
We have
\[
M_{2d}(a_p) = S_{2d}(p)p^{d} + O(p^{d-1}), \qquad
\lim_{p \to \infty} S_{2d}(p) = \frac{(2d)!}{d!(d+1)!}.
\]
\end{theorem}
Note that $\frac{(2d)!}{d!(d+1)!}$ is an integer, namely the $d$-th Catalan number.
It is worth noting that this sort of geometric averaging can be generalized rather substantially; see \cite{katz}.

\subsection{Fixed curve, varying prime}

We now turn to our primary question of interest. Let us now take $E$ to be a fixed elliptic curve\footnote{By definition, an elliptic curve over $K$ is a curve of genus 1 equipped with the choice of a $K$-rational point. For a general curve of genus 1, the existence of a rational point is automatic when $K$ is a finite field but not when $K$ is a number field; in the latter case, though, the Jacobian construction produces a true elliptic curve with the same zeta functions. Consequently, there is no real gain in generality by considering the Sato-Tate conjecture for genus 1 curves over $K$.} over a number field $K$. For $\frakq$ a prime ideal of $K$, we write $q$ for the absolute norm of $K$ and $\FF_q$ for the residue field of $\frakq$. We will only be interested in primes at which $\frakq$ has \emph{good reduction}, so that one may reduce modulo $\frakq$ to get an elliptic curve $E_q$ over $\FF_q$. (This excludes only finitely many primes; for instance, if $E$ is given by a Weierstrass equation, then each prime ideal not dividing $\Delta$ is a prime of good reduction.) 

We again write $\#E(\FF_q)$ in the form $q+1-a_{\frakq}$. (We write $a_{\frakq}$ instead of $a_q$ because there are in general several prime ideals of the same norm, which will have different point counts.)
Since $a_{\frakq}$ now ranges over an interval of varying length, it is natural to renormalize it by defining the new quantity
\[
\overline{a}_{\frakq} = \frac{a_{\frakq}}{\sqrt{q}} \in [-2,2].
\]
One can then ask about the distribution of $\overline{a}_{\frakq}$ as $\frakq$ varies ``randomly,'' but it is tricky to define probability measures on infinite sets. Instead, we should consider all the prime ideals with $q \leq N$ as a probability space with the uniform distribution, then ``take the limit as $N \to \infty$.'' This last step still needs some precision,
but as a start we can make histogram plots for various values of $N$, then view the animated picture given by varying $N$. For some examples of such a visualization, see
\begin{center}
\url{http://math.mit.edu/~drew/g1SatoTateDistributions.html}
\end{center}
Two more precise constructions are the following.
\begin{itemize}
\item
For each positive integer $d$, consider the moment $M_d(\overline{a}_{\frakq}; N)$ over the set of primes up to $N$ as a function of $N$, then take the limit as $N \to \infty$ to get a sequence of ``limiting moments''
$M_d(\overline{a}_{\frakq})$.
\item
Produce a measure $d\mu$ on the interval $[-2,2]$ with the following \emph{equidistribution property}: for any continuous function $f: [-2,2] \to \RR$, we have\footnote{We write $\sum_{q \leq N}$ as shorthand for running over the prime ideals $\frakq$ of $K$ of norm at most $N$. Since there are only finitely many primes of bad reduction, we may define $\overline{a}_{\frakq}$ for them arbitrarily without changing the limit.}
\[
\lim_{N \to \infty} \frac{\sum_{q \leq N} f(\overline{a}_{\frakq})}{\sum_{q \leq N} 1}
= \int_{-2}^2 f\,d\mu.
\]
This can be imagined as a ``time average equals space average'' property in the spirit of ergodic theory.
\end{itemize}
These points of view are equivalent in the following sense: under some technical hypotheses which we omit (but are mild enough not to be a a concern in our applications), 
the existence of the limiting moments implies equidistribution for some $\mu$ and vice versa, and moreover
\[
M_d(\overline{a}_{\frakq}) = \int_{-2}^2 \overline{a}_{\frakq}^d \,d\mu;
\]
that is, the $d$-th limiting moment is in fact the $d$-th moment of the limiting measure.

Using either point of view, if one experiments with various elliptic curves over $K$, one
observes exactly three possible limiting behaviors. The key distinguishing feature here is complex multiplication: the endomorphism ring of $E$, viewed as an algebraic group over $\overline{K}$ (or a complex torus), is either the trivial ring $\ZZ$ (which is present for any commutative algebraic group) or an order in an imaginary quadratic field $M$. In the latter case, we say that $E$ has \emph{complex multiplication in $M$}.

We start with the case of complex multiplication in a subfield of $K$. In this case, 
everything may be deduced easily from Hecke's description of the $a_{\frakq}$ in terms of Grossencharacters. See for instance \cite[Chapter~1]{silverman2}.

\begin{theorem}
Suppose that $E$ has complex multiplication in $M$ and that $M \subseteq K$.
Then the limiting moments exist and satisfy 
\[
M_{2d}(\overline{a}_{\frakq}) = \frac{(2d)!}{(d!)^2},
\]
and equidistribution occurs for the continuous measure defined by the function
\[
\frac{1}{\pi} \frac{1}{\sqrt{4-\overline{a}_{\frakq}^2}}
\]
(which is not bounded at the endpoints, but the improper integral converges).
\end{theorem}

When the complex multiplication occurs in a field not contained in $K$ (as for instance is always the case when $K = \QQ$), one gets a similar but slightly different statement, with essentially the same proof.
\begin{theorem}
Suppose that $E$ has complex multiplication in $M$ and that $M \not\subseteq K$.
Then the limiting moments exist and satisfy 
\[
M_{2d}(\overline{a}_{\frakq}) = \frac{(2d)!}{2(d!)^2},
\]
and equidistribution occurs for the average of the continuous measure defined by
\[
\frac{1}{\pi} \frac{1}{\sqrt{4-\overline{a}_{\frakq}^2}}
\]
with a Dirac point measure concentrated at $\overline{a}_{\frakq} = 0$.
\end{theorem}
The presence of the Dirac component, and the fact that it contributes half of the total measure, can be explained by the fact that in this setting $a_\frakq = 0$ exactly when $\frakq$ is a prime which remains inert in the compositum $KM$.

We now turn to the case where $E$ does not have complex multiplication, which is a certain sense is the most typical: a ``randomly chosen'' elliptic curve does not have complex multiplication. In this case, numerical evidence (first collected by Sato for $K = \QQ$)
and theoretical considerations (first raised by Tate, later expanded by Serre; see below)
lead to the following conjecture.
\begin{conjecture}[Sato-Tate]
Suppose that $E$ does not have complex multiplication. Then the limiting moments exist and satisfy
\[
M_{2d}(\overline{a}_{\frakq}) = \frac{(2d)!}{d!(d+1)!},
\]
and equidistribution occurs for the continuous measure defined by the function
\[
\frac{4}{\pi} \sqrt{4 - \overline{a}_{\frakq}^2}.
\]
\end{conjecture}
Note the recurrence of the Catalan numbers; this means that the limiting distribution in this case coincides with the limiting distribution computed by Birch, in which one first considers all possible elliptic curves over $\FF_p$ for a fixed prime $p$, then takes the limit as $p \to \infty$.

The great difficulty with this conjecture is that one does not have a description of the $a_q$ as simple as that given by Hecke in the case of complex multiplication.
In the case $K = \QQ$, one does at least have the \emph{modularity of elliptic curves} (proved by Wiles in the case of everywhere semistable reduction \cite{wiles, taylor-wiles} and Breuil--Conrad--Diamond--Taylor in general \cite{bcdt}), which asserts that the $a_p$ arise as Fourier coefficients of a certain modular form.
Although this by itself is not enough to establish the conjecture, with a great deal of additional work (more on which below) one can prove the following.
\begin{theorem}
The Sato-Tate conjecture holds whenever the field $K$ is totally real.
\end{theorem}

\subsection{Interlude: the Chebotarev density theorem}

In order to generalize the Sato-Tate conjecture to other abelian varieties, especially the Jacobians of genus 2 curves, we will need to take an alternate point of view. This point of view also happens to provide an indication of how to approach the proof of equidistribution.

We begin by recalling another equidistribution property in number theory: the Chebotarev density theorem. Let $f(T)$ be an irreducible (nonconstant) polynomial of degree $n$ over a number field $K$. For each prime ideal $\frakq$ of $K$ aside from finitely many exceptions, the reduction of $f$ mod $\frakq$ will still have degree $n$ and will factor into distinct irreducible subfactors.
If we read off the degrees of these factors, we obtain an unordered partition $\pi_\frakq$ of $n$.
\begin{theorem}
Let $L$ be the splitting field of $f(T)$ over $K$, put $G = \Gal(L/K)$, and consider the permutation action of $G$ on the roots of $f$ in $L$.
For each unordered partition $\pi$ of $n$, let
 $c_\pi$ be the probability that a random element of $G$ acts with cycle structure $\pi$. Then
\[
\lim_{N \to \infty} \frac{\sum_{q \leq N, \pi_\frakq = \pi} 1}{\sum_{q \leq N} 1}
= c_\pi.
\]
\end{theorem}

To see how this follows from the usual Chebotarev theorem, we may first view cycle structure as a function on $G$ which factors through the quotient $h: G \to \Conj(G)$ to the space of conjugacy classes of $G$. Let us equip $G$ with the uniform measure and $\Conj(G)$ with the \emph{uniform measure}, so that  
for any function $g: \Conj(G) \to \RR$ we have
\[
\int_{\Conj(G)} g = \int_G (g \circ h).
\]
In concrete terms, each class in $\Conj(G)$ is weighted proportionally to its cardinality.

For each prime ideal $\frakq$ of $K$ with finitely many exceptions (namely the primes that ramify in the splitting field of $f$), we may define a \emph{Frobenius class} $g_{\frakq} \in \Conj(G)$ 
using Artin's construction: for any prime ideal $\tilde{\frakq}$ of $L$ lying over $\frakq$,
there is a unique element $g \in G$ stabilizing $\tilde{\frakq}$ and satisfying
\[
x^{g} \equiv x^q \pmod{\tilde{\frakq}} \qquad (x \in \frako_L),
\]
and $g_{\frakq}$ is the conjugacy class of $g$. We then have the usual statement of the density theorem.
\begin{theorem}[Chebotarev]
The elements $g_{\frakq}$ are equidistributed in $\Conj(G)$ for the image of the uniform measure on $G$.
\end{theorem}
We say nothing about the proof except to point out that it can be interpreted as an example of Serre's equidistribution formalism, for which see the end of the second lecture.

\subsection{A suggestive reformulation}

In the previous discussion, we resolved a problem about the distribution of factorizations of reductions of a fixed polynomial modulo a varying prime by ``lifting'' this problem to a problem about the distribution of a sequence of conjugacy classes within a finite group. We now make a first attempt to perform such a lifting for the problem of equidistribution of the values of $\overline{a}_{\frakq}$ arising from an elliptic curve over $K$. This attempt will remain somewhat incomplete until we give a formal definition of Sato-Tate groups in the second lecture.

The condition that $\overline{a}_{\frakq} \in [-2,2]$ is equivalent to asserting that the polynomial $\overline{P}(T) = T^2 - \overline{a}_\frakq T+ 1$ has roots which are complex conjugates and lie on the unit circle. It is also equivalent to asserting that $\overline{P}(T)$ is the characteristic polynomial of some matrix in the Lie group $\SU(2)$ of unitary matrices with determinant 1. Recall that explicitly,
\[
\SU(2) = \left\{ A
\in \GL_2(\CC): A^{-1} = A^*, \det(A) = 1 \right\}
\] 
where $A^* = \overline{A}^T$ is the conjugate transpose of $A$.
In this case, the matrices in $\SU(2)$ with characteristic polynomial $\overline{P}(T)$ form a unique conjugacy class; that is, we may identify $\Conj(\SU(2))$ with $[-2,2]$ via the trace map.

Recall that any compact topological group admits a unique translation-invariant measure called the \emph{Haar measure}; in particular, for a finite group with the discrete topology it is simply the uniform measure.
It is easily checked that under the previous identification, the measure on $[-2,2]$ corresponds to the image measure on $\Conj(\SU(2))$ arising from the Haar measure on $\SU(2)$.  By analogy with Chebotarev, we now see that the Sato-Tate conjecture is equivalent to the equidistribution in $\Conj(\SU(2))$ of the conjugacy classes with traces $\overline{a}_{q}$!

What about the cases of complex multiplication? In the case where $M \subseteq K$, the limiting measure can be interpreted as the image via the trace map of the Haar measure not on $\SU(2)$, but rather on the subgroup $\SO(2)$ of rotations. 
In the case where $M \not\subseteq K$, one instead gets the image of the Haar measure on the normalizer of $\SO(2)$ in $\SU(2)$, which is no longer a connected group: it has two connected components, one consisting of $\SO(2)$ itself, and the other on which the trace is identically zero. This suggests that we should be able to naturally ``lift'' the equidistribution problems in these two cases to problems about the equidistribution of certain conjugacy classes in these two groups; we will see how to do this in the second lecture.

\section{The Sato-Tate group of an abelian variety}

In this lecture, we describe a generalization of the Sato-Tate conjecture to arbitrary abelian varieties. The key point is to identify a compact Lie group and a sequence of conjugacy classes corresponding to the prime ideals over the base field, and then the conjecture simply becomes the equidistribution property. We also include some discussion of a proof strategy for the Sato-Tate conjecture and the extent to which it succeeds.

\subsection{The zeta function of an abelian variety over a finite field}

In order to even ask the correct question in the case of abelian varieties, we must first recall Weil's definition of, and results about, the zeta function of an abelian variety over a finite field. See for example \cite{freitag-kiehl} for a complete presentation.

For $X$ a variety over a finite field $\FF_q$, the  \emph{zeta function} of $X$ is the
complex-analytic function defined in a suitable right half-plane by the absolutely convergent product
\[
\zeta(X,s) = \prod_x (1 - q^{-s \deg(x)})^{-1},
\]
where $x$ runs over the closed points of $X$ (equivalently, the Galois orbits of $\overline{\FF}_q$-points) and $\deg(x)$ denotes the degree of $x$ (i.e., the cardinality of a Galois orbit). It is often more useful to compute $\zeta(X,s)$ via the identity
\begin{equation} \label{eq:zeta identity}
\zeta(X,s) = \exp \left( \sum_{n=1}^\infty \frac{q^{-ns}}{n} \#X(\FF_{q^n}) \right)
\end{equation}
of formal power series in $q^{-s}$.

The general properties of $\zeta(X,s)$ were established via the collective efforts of Dwork, Grothendieck, Deligne, et al. during the 1960's. However, we will only be interested in two special cases established by Weil himself somewhat earlier (which provided motivation for the more general results).

\begin{theorem}[Weil] \label{T:weil1}
Let $C$ be a (smooth, projective, geometrically connected) curve over $\FF_q$ of genus $g$. Then we have
\[
\zeta(C,s) = \frac{P(q^{-s})}{(1-q^{-s})(1-q^{1-s})}
\]
for some polynomial $P(T) = P_C(T) \in \ZZ[T]$ with the following properties.
\begin{itemize}
\item
We have $P(0) = 1$ and $\deg(P) = 2g$.
\item
For $\overline{P}(T) = P(T/\sqrt{q})$, we have the identity
\[
\overline{P}(1/T) = T^{-2g} \overline{P}(T).
\]
\item
If we factor $\overline{P}(T)$ in $\CC$ as $(1 - \alpha_1 T)\cdots (1 - \alpha_{2g}T)$,
then $\left| \alpha_i \right| = 1$ for $i=1,\dots,2g$.
\end{itemize}
\end{theorem}
Some remarks about this statement:
\begin{itemize}
\item
Thanks to \eqref{eq:zeta identity},
\[
\#C(\FF_{q^n}) = q^n+1 - q^{n/2}(\alpha_1^n + \cdots + \alpha_{2g}^n).
\]
\item 
If $C$ is an elliptic curve, then $g=1$, so \eqref{eq:zeta identity} implies
\[
P(T) = 1 - a_q T + q T^2.
\]
In particular, Weil's theorem implies Hasse's theorem.
\item
The complex zeroes of $\zeta(C,s)$ lie on the line $\mathrm{Re}(s) = \frac{1}{2}$.
Thus Weil's theorem provides a bridge between Hasse's theorem and the Riemann hypothesis. \end{itemize}

\begin{theorem}[Weil] \label{T:weil2}
Let $A$ be an abelian variety over $\FF_q$ of dimension $g$. Then we have
\[
\zeta(X,s) = \frac{P_1(q^{-s})\cdots P_{2g-1}(q^{-s})}{P_0(q^{-s}) \cdots P_{2g}(q^{-s})}
\]
where $P(T) = P_A(T)$ satisfies the same conditions as in Theorem~\ref{T:weil1}
and
\[
P_k(T) = \prod_{1 \leq i_1 < \cdots < i_k \leq 2g} (1 - q^{k/2} \alpha_{i_1} \cdots \alpha_{i_k} T).
\]
Moreover, if $A$ is the Jacobian variety associated to a curve $C$, then $P_A(T) = P_C(T)$.
\end{theorem}
Thanks to \eqref{eq:zeta identity},
\[
\#A(\FF_{q^n}) = (1 - q^{n/2} \alpha_1^n) \cdots (1 - q^{n/2} \alpha_{2g}^n).
\]

\subsection{Computing zeta functions}

The analogue of the Sato-Tate conjecture for an abelian variety over a number field will involve the distribution of not just the number of points, but the full zeta function of the reduction modulo various prime ideals. Of course, in order to make reasonable conjectures, one needs to have access to numerical data; although the collection of such data is not the subject of this lecture series, at least a few remarks are in order.

Let us continue to take $A$ to be an abelian variety over a finite field $\FF_q$.
By Theorem~\ref{T:weil2}, the full zeta function of $A$ is determined by the 
polynomial $P_A(T)$. Using \eqref{eq:zeta identity} and the symmetry property of $\overline{P}$, one can compute $P_A(T)$ given $\#A(\FF_{q^n})$ for $n=1,\dots,g$.
The computation of the $\#A(\FF_{q^n})$ is a finite computation in principle (at least assuming $A$ is given in a sufficiently explicit manner), but often quite infeasible in practice.

In the spirit of this conference, let us then restrict attention to the case where $A$ is the Jacobian of a hyperelliptic curve $C$. In this case, a much more practical suite of algorithms can be obtained using $p$-adic analytic techniques. The first of these was introduced by the author in 2001 \cite{kedlaya-hyperelliptic} and will be described in more detail in Balakrishnan's lecture series elsewhere in this volume. This algorithm was originally optimized for finite fields of small characteristic; it was modified to work better in larger characteristics by Harvey \cite{harvey1}. More recently, Harvey \cite{harvey2} discovered a new optimization specific to the case of interest in this paper, where one starts with a hyperelliptic curve over $\QQ$ and considers its reductions modulo all primes up to some bound. This optimization has been put into practice by Harvey--Sutherland with spectacular results \cite{harvey-sutherland}.

\subsection{An equidistribution conjecture in the generic case}

In terms of Weil's zeta functions,
we can now describe the analogue of the Sato-Tate conjecture for a ``generic'' abelian variety over a number field; in the case of elliptic curves, this corresponds to the case of no complex multiplication. The exact origin of this conjecture is a bit unclear, but it
is implicit in Serre's paper on motivic conjectures \cite{serre-motivic}
and explicit in the book of Katz--Sarnak \cite{katz-sarnak}.

Let $A$ be an abelian variety of dimension $g \geq 1$ over a number field $K$.
For each prime ideal $\frakq$ excluding the finitely many primes at which $A$ fails to have good reduction, we may reduce modulo $\frakq$ to obtain an abelian variety $A_q$ over $\FF_q$. Let $P_{\frakq}(T)$ be the polynomial $P_{A_q}(T)$ occurring in the product expression of this zeta function, and define the renormalized polynomial 
\[
\overline{P}_{\frakq}(T) = \overline{P}_{A_q}(T) = P_{\frakq}(T/\sqrt{q}).
\]
By Theorem~\ref{T:weil2}, the roots of $\overline{P}_{\frakq}(T)$ in $\CC$ lie on the unit circle and occur in reciprocal pairs (in particular, each of $+1$ and $-1$ occurs with even multiplicity). These conditions are equivalent to saying that $\overline{P}_{\frakq}(T)$ occurs as the characteristic polynomial of some matrix in the group $\USp(2g)$ of \emph{unitary symplectic} matrices, where explicitly
\[
\USp(2g) = \{A \in \GL_{2g}(\CC): A^{-1} = A^*, \,
A^T J A = J\}
\]
and $J$ is the block diagonal matrix representing the standard symplectic form:
\[
J = \begin{pmatrix} J_1 &  & 0\\
 & \ddots &  \\
 & & J_1 \end{pmatrix},
\qquad
J_1 = \begin{pmatrix} 0 & 1 \\ -1 & 0 \end{pmatrix}.
\]
Moreover, the characteristic polynomial map on $\Conj(\USp(2g))$ is injective, so again $\Conj(\USp(2g))$ may be identified with a space of polynomials. (Note that $\USp(2) = \SU(2)$, so the case $g=1$ of this discussion aligns with our previous discussion of $\SU(2)$.)

\begin{conjecture}[Generalized Sato-Tate conjecture in the generic case]
Suppose that $A$ has ``no extra structure.'' Then the sequence in $\Conj(\USp(2g))$
corresponding to the polynomials $\overline{L}_{\frakq}(T)$ is equidistributed with respect to the image of Haar measure. 
\end{conjecture}

In order to interpret this conjecture, we must clarify what it means for $A$ to have ``no extra structure.'' For $g=1$, this should mean exactly that $A$ does not have complex multiplication. In fact this characterization works for $g \leq 3$: ``no extra structure'' in these cases means exactly that the endomorphism ring of $A$ consists solely of multiplication by integers. We will see how to modify the conjecture in the presence of extra endomorphisms in the next part of this lecture.

For $g \geq 4$, the exclusion of extra endomorphisms is necessary but not sufficient: extra endomorphisms give rise to ``unexpected'' algebraic subvarieties of $A$ of dimension 1, but one must also make similar restrictions on higher-dimensional subravieties. For example, Mumford discovered examples\footnote{As far as we know, no examples of this type are known to be Jacobians of curves of any genus.} in dimension 4 where there are no extra endomorphisms but there are some unexpected subvarieties of $A$ of dimension 2. One can still state a modified Sato-Tate conjecture in such cases, but there are additional technical subtleties which we will only briefly hint at in the next part of this lecture.

\subsection{Construction of the Sato-Tate group}

We now ask the question: if $A$ is an abelian variety over a number field $K$, how do we build a compact Lie group $G$ and a sequence of elements $g_{\frakq} \in \Conj(G)$ having the right equidistribution property whether or not $A$ has any extra structure?
The answer is essentially given by Serre in \cite{serre-motivic}; it is made more explicit by the author and Banaszak in \cite{banaszak-kedlaya}.

To begin with, choose an embedding $K \hookrightarrow \CC$ (this will only affect the answer up to a conjugation). We may then form the base extension $A_{\CC}$ of $A$ from $K$ to $\CC$; since the latter is a complex algebraic variety, it has an associated topological space $A_{\CC}^{\mathrm{an}}$. Let $V = H_1(A^{\mathrm{an}}_{\CC}, \QQ)$ be the rational singular homology of this topological space; it is a $\QQ$-vector space of dimension $2g$.
It also comes equipped with an alternating pairing $\psi$, namely the Riemann form associated to some polarization (the choice of which is immaterial). Let $V_{\CC}$ be the base extension of $V$ from $\QQ$ to $\CC$; then $\psi$ induces a pairing on $V_{\CC}$. We may then choose a symplectic basis of $V_{\CC}$ with respect to $\psi$ and then consider the action of the unitary symplectic group $\USp(2g)$.

When $A$ has extra structure, we need to cut the group $\USp(2g)$ down to something more closely related to the arithmetic of $A$. From the example of $g=1$, we know that the right subgroup will not always be connected; however, as a first approximation let us make a connected subgroup that will be close to the right answer. This will be the largest subgroup $G^\circ$ of $\USp(2g)$ with the following property. For each positive integer $m$, identify the exterior power $\wedge^m V_{\CC}$ with $H_m(A^{\mathrm{an}}_{\CC}, \QQ)$.
Let $W_m$ be the subspace of $\wedge^m V_{\CC}$ spanned by homology classes represented by algebraic subvarieties of $A_{\overline{K}}$ of dimension $m$. Then $G^{\circ}$ is by definition the subgroup of $\USp(2g)$ which fixes $W_m$ for all $m$.
For instance, for $m=2$, this condition means that $G^{\circ}$ commute with the action of the endomorphisms of $A_{\CC}$ on $V_{\CC}$; for $g \leq 3$, this turns out to be sufficient to cut out $G^{\circ}$ (but not for larger $g$; see the discussion above).

Note that the group $G^{\circ}$ is connected and depends only on $A_{\overline{K}}$; it is in fact closely related to the \emph{Mumford-Tate group} of $A$. To get the \emph{Sato-Tate group} $G$ (also called $\ST(A)$), we note that $G_K = \Gal(\overline{K}/K)$ acts on algebraic subvarieties of $A_{\overline{K}}$, and hence on $W_m$. 
We then consider the set of $\gamma \in \USp(2g)$ such that for some $\tau = \tau(\gamma) \in G_K$, for all $m$ the action of $G^\circ$ on $W_m$ coincides with the action of $\tau$.

By construction, the group $G$ has identity connected component $G^\circ$, and the component group $\pi_0(G) = G/G^\circ$ is finite and receives a surjective map from $G_K$. In fact, we may identify $\pi_0(G)$ with $\Gal(L/K)$ where $L$ is the minimal extension of $K$ such that $G_L$ fixes $W_m$ for all $m$. In particular, for $g \leq 3$, $L$ is the minimal field of definition of the endomorphisms of $A_{\overline{K}}$.

\subsection{Construction of conjugacy classes}

Now that we have a candidate group for our desired equidistribution property, it still remains to identify a sequence of conjugacy classes. This requires a nontrivial argument because the map $\Conj(G) \to \Conj(\USp(2g))$ is not injective, so the polynomial $\overline{L}_{\frakq}(T)$ is not enough to determine a class in $\Conj(G)$.

To construct the class associated to a prime ideal $\frakq$, we must switch from singular homology to $\ell$-adic homology, in the form of the Tate module. Let $\ell$ be any prime number. For each positive integer $m$, the group
$A(\overline{K})[\ell^m]$ of $\ell^m$-torsion points of $A$ is isomorphic to $(\ZZ/\ell^m \ZZ)^{2g}$ and carries an action of $G_K$. By taking the inverse limit,
we obtain a $\ZZ_\ell$-module $T_\ell(A)$ which is free of rank $2g$ and admits a continuous action of $G_K$.

For each $m$, we obtain a Frobenius conjugacy class in $\Aut(A(\overline{K})[\ell^m])$; taking inverse limits gives rise to a class $\tilde{g}_{\frakq} \in \Conj(\Aut(T_\ell(A)))$. Note that this class respects the Weil pairing, so it defines a symplectic automorphism of $T_{\ell}(A)$.

It remains to move this $\ell$-adic construction over to $\CC$ to obtain a class in $\ST(A)$; for this, some trickery is unavoidable. We start by choosing an algebraic (but in no way continuous) embedding of the field $\QQ_{\ell}$ into $\CC$. This allows to map $\tilde{g}_{\frakq}$ to a conjugacy class in $\Sp(V_{\CC}, \psi)$; we may then semisimplify and then divide by the scalar $\sqrt{q}$ to obtain a new class $\overline{g}_{\frakq}$. 
This gives a conjugacy class in a maximal compact subgroup of the Zariski closure of $\ST(A)$ (called the \emph{algebraic Sato-Tate group} in \cite{banaszak-kedlaya});
using properties\footnote{The key properties: any semisimple element with eigenvalues on the unit circle belongs to a maximal compact subgroup; any two maximal compact subgroups are conjugate; any two elements of a single maximal compact subgroup which become conjugate in the ambient group are already conjugate within the compact subgroup.} of maximal compact subgroups of a reductive algebraic group over $\CC$
we obtain a well-defined class $g_{\frakq}$ in $\ST(A)$ itself.

\begin{conjecture}[Generalized Sato-Tate conjecture] \label{T:genST}
For any abelian variety $A$ over the number field $K$ (and any prime number $\ell$), the sequence $g_{\frakq}$ defined above is equidistributed in $\Conj(\ST(A))$ for the image of Haar measure.
\end{conjecture}
As stated, this conjecture includes the \emph{Mumford-Tate conjecture} that the image of 
$G_K$ on $\Aut(T_\ell(A))$ is of finite index in the largest possible subgroup consistent with geometric constraints. This weaker conjecture is substantially easier than the generalized Sato-Tate conjecture; for instance, for elliptic curves it was proved by Serre \cite{serre-abelian} in an even stronger form (considering the action on all of the $T_{\ell}$ at once to obtain a comparison of adelic groups) long before any progress was made on the Sato-Tate conjecture. However, while this weaker conjecture is known
in many more cases than the generalized Sato-Tate conjecture, it is open in the most general case.

\subsection{A word on moments}

In the case $g=1$, we were able to use moments to distinguish among distributions of $\overline{a}_{\frakq}$. For general $g$, we are looking at a $g$-dimensional space of polynomials, so we cannot hope to distinguish these using statistics based on a single function.

Write
\[
\overline{P}_{\frakq}(T) = 1 + \overline{a}_{\frakq,1} T + \cdots + T^{2g}.
\]
One can reconstruct the entire joint distribution of
$\overline{a}_{\frakq,1},\dots, \overline{a}_{\frakq,g}$
from the \emph{joint moments}
\[
\mathbb{E}(\overline{a}_{\frakq,1}^{d_1} \cdots \overline{a}_{\frakq,g}^{d_g})
\qquad (d_1,\dots,d_g = 0,1,\dots),
\]
but not from the moments of the $\overline{a}_{\frakq,i}$ individually. 
That being said, if one already has a classification in hand and one needs only to distinguish the distributions at hand (e.g., if one is seeking to match an abelian surface against the classification of \cite{fkrs} based on numerical data), many fewer moments are needed. It may also be useful to consider additional data; for instance, the group-theoretic distributions we will encounter always split as sums of discrete measures and continuous measures, and one may wish to keep track of the discrete components separately: in principle these are determined by the moment sequence, but one may need some high moments for this and these are difficult to compute numerically.

\subsection{A proof strategy for equidistribution in groups}

We now have many examples of problems of the following form: given a compact Lie group $G$ and a sequence $g_{\frakq}$ of elements of $\Conj(G)$, prove that this sequence is equidistributed for the image of Haar measure. It turns out that such problems can be solved using a general argument analogous to Dirichlet's proof of the uniform distribution of primes across residue classes (and with the proofs of the prime number theorem by Hadamard and de la Vall\'ee Poussin). This argument was probably part of Tate's motivation for making the Sato-Tate conjecture, but was first written out explicitly by Serre
\cite[Chapter~I, Appendix]{serre-abelian}. 

For each finite-dimensional linear representation $\rho$ of $G$, define a function $L(\rho,s)$ on the complex half-plane $\mathrm{Re}(s) > 1$ via the absolutely convergent product
\[
L(\rho, s) = \prod_\frakq \det(1 - \rho(\tilde{g}_\frakq)q^{-s})^{-1},
\]
where $\tilde{g}_{\frakq} \in G$ represents any lift of $g_{\frakq} \in \Conj(G)$ (this choice does not affect the definition).
For $\rho$ the trivial representation, $L(\rho,s)$ is essentially the Dedekind zeta function of $K$ (modulo the omission of finitely many factors), and so it has meromorphic continuation to $\CC$ with a simple pole at $s=1$. Using standard representation theory for compact Lie groups (especially the Peter-Weyl theorem), Serre checks the following.
\begin{theorem}[Serre]
Suppose that for each each nontrivial irreducible representation $\rho$, the function $L(\rho,s)$ admits an analytic continuation to a neighborhood of the point $s=1$ which is (holomorphic and) nonzero at $s=1$.
Then the elements $g_{\frakq}$ are equidistributed in $\Conj(G)$ for the image of the Haar measure on $G$.
\end{theorem}

As an aside, we note that when $\rho$ is not irreducible, $L(\rho,s)$ has a pole at $s=1$ of order equal to the multiplicity of the trivial representation in $\rho$. This fact can be used to give a representation-theoretic computation of some moment sequences, such as for the trace in a fixed representation. (In particular, this explains why the moment sequences we have seen previously consist of nonnegative integers.)

Let us briefly discuss how this formalism applies in the cases we have mentioned previously.
\begin{itemize}
\item
For the Chebotarev density theorem, we may check Serre's criterion using the known analytic properties of Dirichlet $L$-functions (plus class field theory to guarantee that this accounts for all $L$-functions associated to abelian Galois groups) and Artin's theorem on induced characters.
\item
For elliptic curves with complex multiplication, we may check Serre's criterion (modulo the exact definition of the equidistribution problem, for which see the second lecture) using Hecke's description of the $a_{\frakq}$ in terms of Grossencharacters, whose $L$-functions are well understood by analogy with the classical theory of the Riemann zeta function.
\item
For elliptic curves without complex multiplication, the best strategy currently known for checking Serre's criterion is to use potential modularity results of Taylor \cite{taylor-potential} for the Galois representations associated to the irreducible representations of $\SU(2)$, namely the symmetric powers of the Tate module. These results represent a significant advance over the original modularity theorem for elliptic curves, which could only handle the first symmetric power. For $K = \QQ$, this strategy was executed by Harris--Shepherd-Barron--Taylor and Clozel--Harris--Taylor \cite{hsbt, cht}; for $K$ totally real, it has been carried out by
Barnet-Lamb--Geraghty--Gee \cite{blgg}.

\item
For abelian varieties of dimension $g>1$, it appears unlikely that current technology will suffice to prove the generalized Sato-Tate conjecture in any cases without extra structure.
However, one can still hope to treat cases with sufficient extra structure. For example,
for abelian varieties with complex multiplication, one can again use Hecke's arguments; this has been worked out explicitly\footnote{Note that \cite{johansson} actually proves equidistribution using a different definition of the Sato-Tate group, then verifies \emph{a posteriori} Serre's criterion for the Sato-Tate group as we have defined it. In fact the two constructions are expected to coincide, 
but this is not checked in \cite{johansson}.} by Johansson \cite{johansson}.
One can also extend the methods used for elliptic curves to certain classes of abelian surfaces; see again \cite{johansson}.
\end{itemize}

It is worth noting that in the presence of more detailed information about the zeroes and poles of $L(\rho, s)$, one can improve the equidistribution statement by controlling the error term. See for example \cite{bucur-kedlaya}.

\section{Sato-Tate groups of abelian surfaces}

In this lecture, we specialize to the case of abelian surfaces and describe the  classification of Sato-Tate groups of abelian surfaces given in \cite{fkrs}.
This can be viewed as giving an explicit form to the generalized Sato-Tate conjecture 
(Theorem~\ref{T:genST}) in the case $g=2$; however, using the preceding definition of the Sato-Tate group, the classification theorem is in fact unconditional.

Note that as for $g=1$, Andrew Sutherland has produced visualizations of equidistribution
for various genus 2 curves. See
\begin{center}
\url{http://math.mit.edu/~drew/g2SatoTateDistributions.html}
\end{center}

Throughout this lecture, let $A$ be an abelian surface over a number field $K$.
One may wish to keep in mind the case of the Jacobian of a genus 2 (necessarily hyperelliptic) curve $C$ over $K$, but most of our work will be on $A$ rather than $C$.

\subsection{Overview of the classification: connected parts}

We start with the well-known classification of Mumford-Tate groups of abelian surfaces, rephrased in terms of the connected part of the Sato-Tate group.
\begin{theorem}
Up to conjugation within $\USp(4)$, there are exactly $6$ groups that occur as the connected parts of Sato-Tate groups of abelian surfaces over number fields:
\[
\SO(2), \SU(2), \SO(2) \times \SO(2), \SO(2) \times \SU(2), \SU(2) \times \SU(2), 
\USp(4).
\]
Moreover, all $6$ groups occur for abelian surfaces over $\QQ$.
\end{theorem}

To see that all of these options occur, let $E_1,E'_1$ be nonisogenous elliptic curves over $K$ with complex multiplication; let $E_2,E'_2$ be nonisogenous elliptic curves over $K$ without complex multiplication; and let $A$ be an abelian surface such that $A_{\overline{K}}$ has trivial endomorphism ring. (For an explicit example of $A$, we may take the Jacobian of the curve $y^2 = x^5 - x + 1$, as originally shown by Zarhin \cite{zarhin}.)
Then the connected parts of the Sato-Tate groups of the abelian surfaces
\[
E_1 \times E_1, E_2 \times E_2, E_1 \times E'_1, E_1 \times E_2, E_2 \times E'_2, A
\]
are exactly the ones listed in the theorem. However, it is also possible to realize all of the connected parts using absolutely simple abelian varieties; we will explain this later in the lecture.

\subsection{Overview of the classification: component groups}

We next discuss the component groups of Sato-Tate groups of abelian surfaces. We start with one form of the main result of \cite{fkrs}.

\begin{theorem}[Fit\'e-Rotger-Kedlaya-Sutherland]
Up to conjugation within $\USp(4)$, there are exactly $52$ groups that occur as Sato-Tate groups of abelian surfaces over number fields, all of which can be realized using genus 2 curves. Of these groups, exactly $34$ groups occur for abelian surfaces over $\QQ$, all of which can be realized using genus 2 curves over $\QQ$. (There is a $35$-th group that can occur over some totally real fields, but not over $\QQ$.)
\end{theorem}

To give more information about these groups, we may consider each of the 6 options for the connected part and then list the options for the component group. See Table~\ref{table:component groups}.
\begin{table}[ht] \label{table:component groups}
\begin{tabular}{c|l}
Connected part & Component groups \\
\hline
$\SO(2)$ & $\Cyc_1, \Cyc_2, \Cyc_2, \Cyc_2, \Cyc_3, \Cyc_4, \Cyc_4, \Cyc_6, \Cyc_6, \Cyc_6$, \\
& $\Dih_2, \Dih_2, \Dih_2, \Dih_3, \Dih_3, \Dih_4, \Dih_4, \Dih_4, \Dih_6, \Dih_6, \Dih_6, \Dih_6, \Alt_4, \Sym_4, \Sym_4$, \\
& $\Cyc_4 \times \Cyc_2, \Cyc_6 \times \Cyc_2, \Dih_2 \times \Cyc_2, \Dih_4 \times \Cyc_2, \Dih_6 \times \Cyc_2, \Alt_4 \times \Cyc_2, \Sym_4 \times \Cyc_2$ \\
$\SU(2)$ & $\Cyc_1, \Cyc_2, \Cyc_2, \Cyc_3, \Cyc_4, \Cyc_6, \Dih_2, \Dih_3, \Dih_4, \Dih_6$\\
$\SO(2) \times \SO(2)$ & $\Cyc_1, \Cyc_2, \Cyc_2, \Cyc_4, \Dih_2$ \\
$\SO(2) \times \SU(2)$ & $\Cyc_1, \Cyc_2$ \\
$\SU(2) \times \SU(2)$ & $\Cyc_1, \Cyc_2$ \\
$\USp(4)$ & $\Cyc_1$
\end{tabular}
\caption{Component groups of Sato-Tate groups. Here $\Cyc_n, \Dih_n, \Alt_n, \Sym_n$ denote the cyclic, dihedral, alternating, and symmetric groups whose standard permutation representations have degree $n$. Multiple listings of the same group indicate distinct extensions within $\USp(2g)$.}
\end{table}

One notices immediately from the table that the number of possible component groups increases as the connected part gets smaller. This makes sense on many levels.
For instance, recall that the component group is identified with $\Gal(L/K)$ for $L$ the minimal field over which the endomorphisms of $A_{\overline{K}}$ can all be realized. A Sato-Tate group with small connected part indicates the presence of many endomorphisms,
which then may fuse together to live over a large extension of $K$. In fact, from the table we see that the largest possible value for $[L:K]$ is $\#(\Sym_4 \times \Cyc_2) = 48$;
an example achieving this bound is the Jacobian of
\[
y^2 = x^6 - 5x^4 + 10x^3 - 5x^2 + 2x - 1.
\]
Previously it was only known that $[L:K]$ always divides $11520 = 2^8 \cdot 3^2 \cdot 5$, by specializing a bound for all $g$ given by Silverberg \cite{silverberg}.

\subsection{Sato-Tate groups and endomorphism algebras}

Let $\End(A_{\overline{K}})$ be the endomorphism ring of $A_{\overline{K}}$. 
\begin{theorem}[Fit\'e-Rotger-Kedlaya-Sutherland]
The following statements hold.
\begin{enumerate}
\item[(a)]
The group $\ST(A)^\circ$ (up to conjugation within $\USp(4)$) uniquely determines, and is uniquely determined by, the $\RR$-algebra $\End(A_{\overline{K}})_{\RR} = \End(A_{\overline{K}}) \otimes_{\ZZ} \RR$.
\item[(b)]
The group $\ST(A)$ (up to conjugation within $\USp(4)$) uniquely determines, and is uniquely determined by, the $\RR$-algebra $\End(A_{\overline{K}})_{\RR}$ equipped with its $G_K$-action.
\end{enumerate}
\end{theorem}

The fact that we must pass from $\End(A_{\overline{K}})$ to $\End(A_{\overline{K}})_{\RR}$
in order to compute the Sato-Tate group means that one can have simple and nonsimple abelian surfaces with the same Sato-Tate group. For example, an absolutely simple abelian surface with complex multiplication in a quartic field has the same connected part of its Sato-Tate group as the product of two nonisogenous elliptic curves with complex multiplication. To derive more examples of this sort, see Table~\ref{table:absolute type}.

\begin{table} \label{table:absolute type}
\begin{tabular}{c|c|c}
Absolute type & $\ST(A)^\circ$ & $\End(A_{\overline{K}})_{\RR}$ \\
\hline
$\mathbf{A}$ & $\USp(4)$ & $\RR$ \\
$\mathbf{B}$ & $\SU(2) \times \SU(2)$ & $\RR \times \RR$ \\
$\mathbf{C}$ & $\SO(2) \times \SU(2)$ & $\RR \times \CC$ \\
$\mathbf{D}$ & $\SO(2) \times \SO(2)$ & $\CC \times \CC$ \\
$\mathbf{E}$ & $\SU(2)$ & $\mathrm{M}_2(\RR)$ \\
$\mathbf{F}$ & $\SO(2)$ & $\mathrm{M}_2(\CC)$
\end{tabular}
\caption{Real endomorphism algebras corresponding to connected parts of Sato-Tate groups. The left column represents the \emph{absolute type} of $A$ in the terminology of \cite{fkrs}.}
\end{table}

In \cite[\S 4]{fkrs}, one finds a detailed enumeration of the possible Galois actions on real endomorphism algebras of abelian surfaces; these are labeled in terms of \emph{Galois type}. Let $L$ be the minimal field of definition of endomorphisms of $\overline{A}_{\overline{K}}$. In most cases, the label of a Galois type has the form
\[
\mathbf{L}[\Gal(L/K)]
\]
where $\mathbf{L} \in \{\mathbf{A},\dots,\mathbf{F}\}$ is the \emph{absolute type} (i.e., the underlying real endomorphism algebra) as labeled in Table~\ref{table:absolute type}). For $\mathbf{L} = \mathbf{D}, \mathbf{E}$, the label $\mathbf{L}[\Cyc_2]$ is ambiguous, so we expand it to 
\[
\mathbf{L}[\Cyc_2, \End(A_{\overline{K}})_{\RR}^{\Cyc_2}].
\]
For $\mathbf{L} = \mathbf{F}$, there are many ambiguous cases, so we disambiguate in a systematic fashion. In this case, the ring $\End(A_{\overline{K}})_{\QQ}$ is a quaternion algebra over some imaginary quadratic field $M$. When $M \subseteq K$, we use labels as above; when $M \not\subseteq K$, we use labels of the form
\[
\mathbf{F}[\Gal(L/K),\Gal(L/KM),\End(A_{\overline{K}})_{\RR}^H].
\]
It is a corollary of the classification that each Galois type receives a unique label under this scheme.

\subsection{Comments on the proof}

The proof of the classification theorem of \cite{fkrs} consists of three main parts.
The first part \cite[\S 3]{fkrs} is a purely group-theoretical classification of subgroups $G$ of $\USp(4)$ (up to conjugation) satisfying certain conditions imposed by Hodge theory (called the \emph{Sato-Tate axioms} in \cite{fkrs}).
\begin{enumerate}
\item[(ST1)] The group $G$ is a closed subgroup of $\USp(4)$.
\item[(ST2)] There exists a homomorphism $\theta: \mathrm{U}(1) \to G^\circ$ such that $\theta(u)$ has eigenvalues $u,u,u^{-1},u^{-1}$ and does not factor\footnote{This second part of the condition is mistakenly omitted from \cite{fkrs}. See \cite[8.2.3.6(i)]{serre-lectures}.} through any closed proper normal subgroup of $G^\circ$.
\item[(ST3)] For each component $H$ of $G$ and each irreducible character $\chi$ of $\GL_4(\CC)$, the average (for Haar measure) of $\chi(\gamma)$ over $H$ is an integer.
\end{enumerate}

This yields a preliminary list of 55 groups which are eligible to occur as Sato-Tate groups. The second part of the proof is a matching argument using the enumeration of Galois types in \cite[\S 4]{fkrs}; this fails to realize 3 of the groups in the original classification, yielding the list of 52 groups. The third part of the group is an enumeration of some genus 2 curves which provably realize all 52 groups \cite[\S 4.8]{fkrs}.

\subsection{Numerical implications}

Let us now take $C$ to be a genus 2 curve over $K$, and suppose that we can compute the polynomials $P_{\frakq}(T)$ for many prime ideals $\frakq$ (this is indeed practical, as discussed at the end of the second lecture). If we write $\overline{P}_{\frakq}(T)$ in the form
\[
T^4 + \overline{a}_{1,\frakq} T^3 + \overline{a}_{2,\frakq} T^2 + \overline{a}_{1,\frakq} T + 1,
\]
then we can compute (approximations of) the limiting moments of both $\overline{a}_{1,\frakq}$ and $\overline{a}_{2,\frakq}$, and according to the generalized Sato-Tate conjecture these should equal the corresponding averages over the Sato-Tate groups. These moment sequences are tabulated in \cite[\S 6]{fkrs} and can be used to distinguish all 52 possible groups. In fact, in many case much less data is necessary; for instance, the group $\USp(4)$ is the only one for which $M_4(\overline{a}_{1,\frakq}) = 3$, the other groups all yielding larger values.

In some cases, it is more useful to supplement the computation of moments with the additional computation of the densities of special values. For $\overline{a}_{1,\frakq}$, only the value 0 can occur with positive density; for $\overline{a}_{2,\frakq}$, the values $-2,-1,0,1,2$ can occur with positive density. These densities are also tabulated in \cite[\S 6]{fkrs}.

It must be pointed out here that this presentation of \cite{fkrs} inverts the order of discovery! The classification began with this kind of numerical investigation; a partial numerical census was made in \cite{kedlaya-sutherland}, but this turned out not to be exhaustive for (at least) two reasons: only the distribution of $\overline{a}_{1,\frakq}$ was considered, which conflates certain distinct Sato-Tate groups; and the search was not exhaustive enough to detect some of the rarest cases.

\subsection{Future directions: beyond dimension $2$}

To conclude, we give some indications of some possible analogues of the classification in \cite{fkrs}.

The most obvious next step would be to consider abelian threefolds. This is now feasible thanks to \cite{harvey-sutherland}, but looks to be a daunting task for two reasons. On one hand, there is a combinatorial explosion in the number of possible cases, especially when the connected part is as small as possible (i.e., cases of complex multiplication). It is likely that the final number of Sato-Tate groups (up to conjugation) will number in the hundreds or thousands. On the other hand, the rarest of these cases may be extremely difficult to detect using a brute-force search, so some intricate analysis is likely to be required to rule them in or out.

It is instead probably more profitable to consider higher-dimensional varieties which have relatively small pieces of cohomology which can be isolated. An example of this can be found in \cite{fks}, where a partial classification is made of motives of weight 3 whose Hodge numbers coincide with the symmetric cube of an elliptic curve. The most generic cases of this form can be achieved using the Dwork pencil of quintic threefolds:
\[
x_0^5 + x_1^5 + x_2^5 + x_3^5 + x_4^5 = \lambda x_0 x_1 x_2 x_3 x_4.
\]

\end{document}